# Reliable Routing of Road-Rail Intermodal Freight Under Uncertainty

Majbah Uddin[1] and Nathan Huynh[2]*

University of South Carolina
Department of Civil and Environmental Engineering
300 Main St, Columbia, SC 29208, USA

Email: [1]muddin@email.sc.edu, [2]nathan.huynh@sc.edu

ORCiD ID: [1]0000-0001-9925-3881, [2]0000-0002-4605-5651

*Corresponding Author Contact Information
Nathan Huynh
University of South Carolina
Department of Civil and Environmental Engineering
300 Main St, Columbia, SC 29208, USA
Telephone: (803) 777-8947
Fax: (803) 777-0670
Email: nathan.huynh@sc.edu

**Abstract**

Transportation infrastructures, particularly those supporting intermodal freight, are vulnerable to natural disasters and man-made disasters that could lead to severe service disruptions. These disruptions can drastically degrade the capacity of a transportation mode and consequently have adverse impacts on intermodal freight transport and freight supply chain. To address service disruption, this paper develops a model to reliably route freight in a road-rail intermodal network. Specifically, the model seeks to provide the optimal route via road segments (highway links), rail segments (rail lines), and intermodal terminals for freight when the network is subject to capacity uncertainties. To ensure reliability, the model plans for reduced network link, node, and intermodal terminal capacity. A major contribution of this work is that a framework is provided to allow decision makers to determine the amount of capacity reduction to consider in planning routes to obtain a user-specified reliability level. The proposed methodology is demonstrated using a real-world intermodal network in the Gulf Coast, Southeastern, and Mid-Atlantic regions of the United States. It is found that the total system cost increases with the level of capacity uncertainty and with increased confidence levels for disruptions at links, nodes, and intermodal terminals.

# 1 Introduction

Freight transportation involves various transportation modes, such as road, rail, air and water. The use of different transportation modes provides greater efficiency because it takes advantages of the strength of each transportation mode. Intermodal freight transportation uses two or more modes to transport goods without handling the goods themselves. Intermodal transportation offers an attractive alternative to unimodal transportation by highway in terms of cost for freight transported over long distances, and it reduces the carbon footprint of transport compared to the highway mode (Bureau of Transportation Statistics 2015). In recent years, intermodal freight transport volume has grown significantly due to the aforementioned advantages.

Transportation infrastructures, particularly those supporting intermodal freight, are vulnerable to natural disasters (e.g., hurricane, earthquake, flooding) and man-made disasters (e.g., accidents, labor strike). These disruptions can drastically degrade the capacity of a transportation mode and consequently have adverse impacts on intermodal freight transport and freight supply chain (Miller-Hooks et al. 2012; Uddin and Huynh 2016). For examples, Hurricane Katrina significantly damaged the transportation infrastructure in the Gulf Coast area (Godoy 2007), and the West Coast port labor strike severely disrupted the U.S. freight supply chain (D'Amico 2002). Therefore, there is a need to develop a modeling framework that takes into account the reliability of the freight transport network when making strategic routing decisions. Network reliability means that the network can continue to deliver acceptable service when faced with disasters or disruptions that reduce capacity of network links, nodes, and intermodal terminals.

The majority of the studies that deal with intermodal freight shipments seek to minimize routing cost. Barnhart and Ratliff (1993) proposed a model for minimizing routing cost in a road-rail intermodal network. They developed procedures involving shortest paths and matching algorithm to help shippers in deciding routing options. Boardman et al. (1997) developed a software-based decision support system to assist shippers to select the best combination of transportation modes considering cost, service level, and the type of commodity. Xiong and Wang (2014) developed a bi-level multi-objective genetic algorithm for the routing of freight with time windows in a multimodal network. Ayar and Yaman (2012) investigated an intermodal multicommodity routing problem where release times and due dates of commodities were pre-scheduled in a planning horizon. Uddin and Huynh (2015) developed a methodology for freight traffic assignment in large-scale road-rail intermodal networks to be used by transportation planners to forecast intermodal freight flows. Rudi et al. (2016) proposed a capacitated multicommodity network flow model for the intermodal freight transportation problem that seeks to minimize transportation costs, carbon emissions, and in-transit holding costs. Their model was validated using industry data from an automotive supplier. Qu et al. (2016) provided a multimode multicommodity service network design model for intermodal freight transportation considering greenhouse gas emissions and intermodal transfer cost.

Intermodal freight transportation has also been studied in the context of network equilibrium. Friesz et al. (1986) presented a network equilibrium model for predicting freight flows considering the role of both shippers and carriers. Their model considered the route choice decisions of shippers and carriers sequentially on a multimodal freight network with nonlinear cost and delay functions. Guelat et al. (1990) used a Gauss–Seidel–Linear approximation algorithm to assign multiproduct freight flows on a multimodal network for strategic planning. In addition to link costs, their algorithm considered intermodal transfer costs to determine the shortest paths. Agrawal and Ziliaskopoulos (2006) used variational inequality to develop a dynamic shipper-carrier freight assignment model. In their model, the market equilibrium is reached when no shipper can reduce its cost by changing carrier for any shipment. Li et al. (2014) proposed a receding horizon approach for the intermodal container flow assignment problem. Their approach assigned container flows using the solution from a nonlinear optimization model that is evaluated at each time step for each node in the network. Corman et al. (2017) presented an equilibrium model for multimodal container transportation. Their model minimized the generalized costs, which is a function of mode, travel time and waiting time for freight consolidation.

In addition to the use of mathematical programs to solve the intermodal freight transportation problem, other methods that have been utilized include network simulation, econometrics, and geospatial analysis. Mahmassani et al. (2007) developed a dynamic freight network simulation-assignment framework that can be used to analyze multiproduct intermodal freight transportation systems. Their framework modeled individual shipment mode-path choice behavior and terminal transfer. Zhang et al. (2008) applied the Mahmassani et al. framework to a pan-European rail network and validated it by analyzing the convergence patterns and performance measures. Lim and Thill (2008) utilized geographic information system-based mapping to evaluate the performance of the U.S. intermodal freight transportation network. They also utilized geographically weighted regressions to identify the factors affecting the improvement of accessibility due to intermodalism. Winebrake et al. (2008) presented a geospatial model for analyzing intermodal freight network in terms of cost, time-of-delivery, energy and environmental impact. Their model can also be used to explore tradeoffs among different mode combinations. Meng and Wang (2010) proposed an algorithm based on Monte Carlo simulation to estimate the probability of shippers selecting an intermodal route involving a port. Their algorithm maximized the utility of shippers, which is defined as a summation of transportation cost and transport time multiplied by the perceived value of time.

All of the aforementioned studies assume that the freight transport network is always functioning and is never disrupted. Daskin (1983) considered disruptions by taking into account the facility unavailability in a maximum covering location problem. Snyder and Daskin (2005) presented a uncapacitated location problem considering failure of facilities in the network. Their reliability models find facility location by taking into account the expected transportation cost after failure, in addition to the minimum operational cost. Cui et al. (2010) extended this work to consider failures with site-dependent probabilities and re-routing of customers when there are failures. Unnikrishnan et al. (2009) developed a two-stage linear program with recourse for the shipper-carrier network design problem under uncertainty. In their model, the shipper decides the optimal capacity for the transshipment nodes in the first stage, and in the second stage, it chooses a routing strategy based on the realized demand. Peng et al. (2011) considered disruptions at facilities in their work on design of reliable logistics network. In contrast, Cappanera and Scaparra (2011) sought to improve network reliability by optimally allocating protective resources in shortest path networks. Chen and Miller-Hooks (2012) developed a method to quantify resilience of an intermodal freight transport network. Miller-Hooks et al. (2012) extended this work to maximize freight transport network resiliency by implementing preparedness and recovery activities within a given budget. Huang and Pang (2014) evaluated resiliency of biofuel transport networks under possible natural disruptions. They formulated a multi-objective stochastic program to optimize the total system cost and total resilience cost. Marufuzzaman et al. (2014) proposed a reliable multimodal transportation network design model, where intermodal hubs are subject to site-dependent disruptions. This model was solved using the accelerated Benders decomposition algorithm and tested on a large-scale network. Uddin and Huynh (2016) proposed a stochastic mixed-integer model for the routing of multicommodity freight in an intermodal network under disruptions. Their study found that goods are better shipped via road-rail intermodal network during disruptions due to the built-in redundancy of the freight transport network. Li et al. (2016) formulated a freight routing model considering both reliability and sustainability under link travel time uncertainty. Fotuhi and Huynh (2017) developed a robust mixed-integer linear model that can be used by railroad operators to evaluate intermodal network expansion options when there are uncertainties in freight demands and network element capacities.

      A number of studies have considered network vulnerability in planning decision. Chen et al. (2007) presented network-based accessibility measures for assessing the vulnerability of transportation networks under disruptions. In addition to the increase in travel time due to link failures, their model considered the behavioral responses of users. Peterson and Church (2008) investigated rail network vulnerability by formulating both uncapacitated and capacitated routing-based model. Garg and Smith (2008) presented a methodology for designing a survivable multicommodity flow network, which analyzes failure scenarios involving multiple arcs. Rios et al. (2000) studied a similar problem, but their objective was to find the minimum-cost capacity-expansion options such that shipments can still be delivered to receivers through the network under disruptions. Gedik et al. (2014) proposed a capacitated

mixed-integer interdiction programming model for coal transportation. They assessed network vulnerability and re-routing of coal by rail under network disruptions. Viljoen and Joubert (2018) presented a model to quantify the impact of transportation infrastructure on supply chain vulnerability; their methodology used a multilayered network.

Another area of research that involves network uncertainty is disaster management, relief routing, response planning, and emergency and humanitarian logistics. Researchers have developed a wide variety of classical optimization programs to address these challenging problems. Haghani and Oh (1996) presented a disaster relief routing model for multicommodity freight in a multimodal network using the concept of time-space network. In the work by Ozdamar et al. (2004), commodity relief routing was studied as a hybrid of classical multicommodity network flow and vehicle routing problem. Given the uncertainty associated with network disruption, their model attempted to deliver commodities such that unsatisfied demand is minimized in a multimodal network. Barbarosoglu and Arda (2004) proposed a stochastic programming model for transporting multicommodity freight through a multimodal network during a natural disaster. Their model considered random arc capacity, where randomness is represented by a finite sample of scenarios. Chang et al. (2007) studied the rescue resources location-routing problem in the event of a flooding disaster. Shen et al. (2009) investigated how to route vehicles in the event of a large-scale bioterrorism emergency. Their solution approach involves adjusting routes generated at the planning level to consider effects of disruptions. Rennemo et al. (2014) proposed a model comprising several stages to optimally locate relief distribution facilities. Cantillo et al. (2018) developed a transportation network vulnerability assessment model to identify critical links for siting distribution centers for disaster response.

**Table 1** Summary of prior studies on the routing of freight

| Study | Mode | Multi commodity | Capacitated link | Delivery deadline | Uncertainty consideration | Probability distribution assumption |
|---|---|---|---|---|---|---|
| Haghani and Oh (1996) | Multiple | ✔ | Multiple | | | |
| Barbarosoglu and Arda (2004) | Road, air | ✔ | Road, air | | ✔ | ✔ |
| Garg and Smith (2008) | Road | ✔ | Road | | ✔ | ✔ |
| Ayar and Yaman (2012) | Road, water | ✔ | Water | ✔ | | |
| Chen and Miller-Hooks (2012) | Road, rail | | Road, rail | ✔ | ✔ | ✔ |
| Miller-Hooks et al. (2012) | Road, rail | | Road, rail | ✔ | ✔ | ✔ |
| Gedik et al. (2014) | Rail | | Rail | | ✔ | ✔ |
| Rudi et al. (2016) | Road, rail, water | ✔ | Road | | | |

| Uddin and Huynh (2016) | Road, rail | ✔ | Road, rail | ✔ | ✔ | ✔ |
| This current study | Road, rail | ✔ | Road, rail | ✔ | ✔ | |

Table 1 provides a summary of the key features addressed by prior studies related to the routing of freight. All of the prior studies where network uncertainty is considered make an explicit assumption about the probability density function (PDF) of the network link and/or node capacity. However, given that disruptive events are rare, there is often limited or no historical data available to determine the PDF of the network link or node capacity under a particular disruption scenario. A wrong assumption could have serious consequences of over design or under design. For example, assuming that a link capacity will follow the normal distribution in the event of a flash flood when in fact it follows a gamma distribution would lead to over design. This study contributes to the current body of knowledge by relaxing this explicit PDF assumption. A novel distribution-free approach is used to provide probabilistic guarantees on the resulting routes. This approach uses symmetric random variation, which is a popular method for solving robust optimization models (Bertsimas and Sim 2004; Ng and Waller 2012).

The remainder of this paper is organized as follows. Section 2 describes the problem considered in this study and presents the mathematical formulation of the problem. Section 3 describes the application of the proposed methodology on a real-world road-rail intermodal transportation network. Lastly, section 4 provides a summary of the study, its contributions, and future work.

## 2 Problem Description and Model Formulation

The main objective of this paper is to develop a reliable routing model for shipment of freight on a road-rail intermodal network that is subject to capacity uncertainty. The problem consists of determining the routes for commodity shipments from their origins (shippers) to destinations (receivers). In this study, it is assumed that the origins and destinations are only accessible via highway links and that every intermodal route will involve at least two intermodal terminals. Additionally, it is assumed that the shipper and receiver facilities are either warehouses or distribution centers and that these facilities do not have rail connections. Figure 1 presents a typical road-rail freight transportation network where shipments can be transported via road-only or intermodal. The network consists of freight shippers, receivers, intermodal terminals, highway links, and rail lines.

Following the notations from Uddin and Huynh (2016), it is assumed that a road-rail intermodal freight transportation network is represented by a directed graph $G = (N, A)$, where $N$ is the set of nodes and $A$ is the set of links. Set $N$ consists of the set of major highway intersections $H$, the set of major rail junctions $R$, and the set of intermodal terminals $S$. Set $A$ consists of the set of highway links $A_h$ and the set of railway links $A_r$. Shipments can change mode at the intermodal terminal nodes, *S*. Each highway link $(i,j) \in A_h$ and railway link $(i,j) \in A_r$ have unit transportation costs associated with them for each commodity $k \in K$ shipment. Each intermodal terminal $s \in S$ has also a unit transfer cost for each commodity $k \in K$ shipment. The definitions of sets, parameters, and decision variables are presented next, followed by the model formulation.

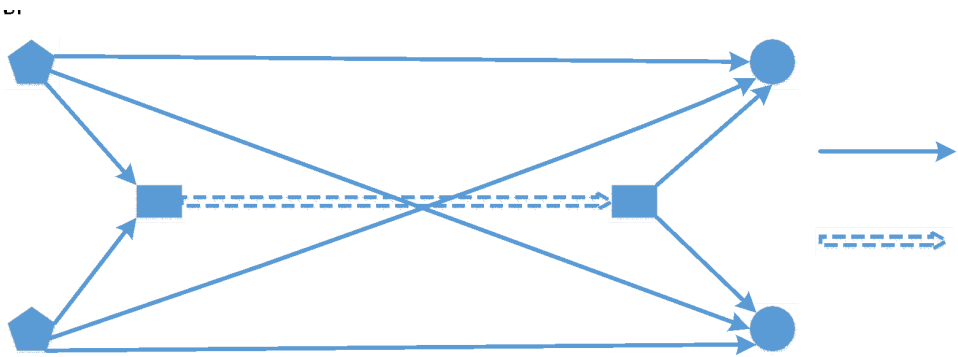

**Fig. 1** An example of road-rail freight transportation network

*Sets/Indices*

| | |
|---|---|
| $H$ | set of major highway intersections |
| $R$ | set of major rail junctions |
| $S$ | set of intermodal terminals |
| $A_h$ | set of highway links |
| $A_r$ | set of railway links |
| $C$ | set of origin-destination (OD) pairs |
| $K$ | set of commodities |
| $P^c$ | set of paths $p$ connecting OD pair $c$ |
| $k$ | commodity type, $k \in K$ |
| $i, j, s$ | node, $i, j, s \in N$ |
| $c$ | an OD pair, $c \in C$ |
| $ori^c$ | origin node of an OD pair $c \in C$ |
| $des^c$ | destination node of an OD pair $c \in C$ |

*Parameters*

| | |
|---|---|
| $d_k^c$ | original demand of commodity $k \in K$ between OD pair $c \in C$ (expressed in number of intermodal containers) |
| $\Psi$ | unit penalty cost for unsatisfied demand |
| $\beta_{ijk}$ | unit cost of transporting commodity $k \in K$ by truck in link $(i,j) \in A_h$ |
| $\tilde{\beta}_{ijk}$ | unit cost of transporting commodity $k \in K$ by rail in link $(i,j) \in A_r$ |
| $\beta_{sk}$ | unit cost of transferring commodity $k \in K$ in intermodal terminal $s \in S$ |
| $Q_{ij}$ | capacity of highway link $(i,j) \in A_h$ |
| $\tilde{Q}_{ij}$ | capacity of railway link $(i,j) \in A_r$ |
| $Q_s$ | capacity of intermodal terminal $s \in S$ |
| $t_{ij}$ | travel time on highway link $(i,j) \in A_h$ |
| $\tilde{t}_{ij}$ | travel time on railway link $(i,j) \in A_r$ |
| $t_s$ | processing time in intermodal terminal $s \in S$ |
| $T_k^c$ | delivery time for commodity $k \in K$ between OD pair $c \in C$ |
| $M$ | sufficiently large number |
| $\varepsilon$ | sufficiently small number |

*Decision variables*

| | |
|---|---|
| $X_{ijk}^c$ | fraction of commodity $k \in K$ transported in highway link $(i,j) \in A_h$ between OD pair $c \in C$ |
| $\tilde{X}_{ijk}^c$ | fraction of commodity $k \in K$ transported in railway link $(i,j) \in A_r$ between OD pair $c \in C$ |

$U_k^c$          unsatisfied demand of commodity $k \in K$ between OD pair $c \in C$

$F_{sk}^c$         fraction of commodity $k \in K$ between OD pair $c \in C$ transferred at terminal $s \in S$

$Y_{sk}^c$         binary variable indicating whether or not intermodal terminal $s \in S$ is selected for commodity $k \in K$ between OD pair $c \in C$ (= 1 if intermodal terminal $s$ is selected for commodity $k$ between OD pair $c$, = 0 otherwise)

$\delta_{ijk}^c$         binary variable indicating whether or not there is any flow in highway link $(i,j) \in A_h$ for commodity $k \in K$ between OD pair $c \in C$ (= 1 if highway link $(i,j)$ carries flow of commodity $k$ between OD pair $c$, = 0 otherwise)

$\tilde{\delta}_{ijk}^c$         binary variable indicating whether or not there is any flow in railway link $(i,j) \in A_r$ for commodity $k \in K$ between OD pair $c \in C$ (= 1 if railway link $(i,j)$ carries flow of commodity $k$ between OD pair $c$, = 0 otherwise)

The multicommodity intermodal freight shipment routing (MIFR) problem is formulated as follows.

$$\text{Minimize} \sum_{c \in C} \sum_{k \in K} \left( d_k^c \left( \sum_{(i,j) \in A_h} \beta_{ijk} X_{ijk}^c + \sum_{(i,j) \in A_r} \tilde{\beta}_{ijk} \tilde{X}_{ijk}^c + \sum_{s \in S} \beta_{sk} F_{sk}^c \right) + \Psi U_k^c \right) \quad (1)$$

subject to

$$\sum_{(i,m) \in A_h} X_{imk}^c - \sum_{(m,i) \in A_h} X_{mik}^c \begin{cases} \leq +1 & \text{if } i = \text{ori}^c \\ \geq -1 & \text{if } i = \text{des}^c, \\ = 0 & \text{otherwise} \end{cases} \forall i \in H, k \in K, c \in C \quad (2)$$

$$\sum_{(i,m) \in A_h : i = \text{ori}^c} X_{imk}^c - \sum_{(m,j) \in A_h : j = \text{des}^c} X_{mjk}^c = 0, \ \forall k \in K, c \in C \quad (3)$$

$$X_{imk}^c \leq \delta_{imk}^c, \ \forall (i,m) \in A_h, k \in K, c \in C \quad (4)$$

$$X_{mik}^c + \delta_{imk}^c \leq 1, \ \forall (m,i) \in A_h, k \in K, c \in C \quad (5)$$

$$\sum_{(i,m)\in A_h} X^c_{imk} \geq M \sum_{(m,i)\in A_h} X^c_{mik}, \quad \forall i \in \text{ori}^c, k \in K, c \in C \tag{6}$$

$$\sum_{(i,n)\in A_r} \widetilde{X}^c_{ink} - \sum_{(n,i)\in A_r} \widetilde{X}^c_{nik} = 0, \quad \forall i \in R, k \in K, c \in C \tag{7}$$

$$\sum_{(s,m)\in A_h} X^c_{smk} - \sum_{(m,s)\in A_h} X^c_{msk} + \sum_{(s,n)\in A_r} \widetilde{X}^c_{snk} - \sum_{(n,s)\in A_r} \widetilde{X}^c_{nsk} = 0, \quad \forall s \in S, k \in K, c \in C \tag{8}$$

$$-M\, Y^c_{sk} \leq \sum_{(s,n)\in A_r} \widetilde{X}^c_{snk} - \sum_{(n,s)\in A_r} \widetilde{X}^c_{nsk} \leq M\, Y^c_{sk}, \quad \forall s \in S, k \in K, c \in C \tag{9}$$

$$-F^c_{sk} \leq \sum_{(s,m)\in A_h} X^c_{smk} - \sum_{(m,s)\in A_h} X^c_{msk} \leq F^c_{sk}, \quad \forall s \in S, k \in K, c \in C \tag{10}$$

$$\sum_{(i,j)\in (A_h \cap p)} \delta^c_{ijk} t_{ij} + \sum_{(i,j)\in (A_r \cap p)} \widetilde{\delta}^c_{ijk} \widetilde{t}_{ij} + \sum_{s\in (S \cap p)} Y^c_{sk} t_s \leq T^c_k, \quad \forall p \in P^c, k \in K, c \in C \tag{11}$$

$$\sum_{c\in C}\sum_{k\in K} d^c_k X^c_{ijk} \leq Q_{ij}, \quad \forall (i,j) \in A_h \tag{12}$$

$$\sum_{c\in C}\sum_{k\in K} d^c_k \widetilde{X}^c_{ijk} \leq \widetilde{Q}_{ij}, \quad \forall (i,j) \in A_r \tag{13}$$

$$\sum_{c\in C}\sum_{k\in K} d^c_k F^c_{sk} \leq Q_s, \quad \forall s \in S \tag{14}$$

$$d^c_k \left(1 - \sum_{i\in H} X^c_{ijk}\right) = U^c_k, \quad \forall k \in K, c \in C, j = \text{des}^c \tag{15}$$

$$\varepsilon \delta^c_{ijk} \leq X^c_{ijk} \leq \delta^c_{ijk}, \quad \forall (i,j) \in A_h, k \in K, c \in C \tag{16}$$

$$\varepsilon \widetilde{\delta}^c_{ijk} \leq \widetilde{X}^c_{ijk} \leq \widetilde{\delta}^c_{ijk}, \quad \forall (i,j) \in A_r, k \in K, c \in C \tag{17}$$

$$\varepsilon Y^c_{sk} \leq F^c_{sk} \leq Y^c_{sk}, \quad \forall s \in S, k \in K, c \in C \tag{18}$$

$$X^c_{ijk} \in [0,1], \quad \forall (i,j) \in A_h, k \in K, c \in C \tag{19}$$

$$\widetilde{X}^c_{ijk} \in [0,1], \quad \forall (i,j) \in A_r, k \in K, c \in C \tag{20}$$

$$F_{sk}^c \in [0,1], \quad \forall s \in S, k \in K, c \in C \tag{21}$$

$$U_k^c \in \mathbb{Z}^+, \quad \forall k \in K, c \in C \tag{22}$$

$$\delta_{ijk}^c \in \{0,1\}, \quad \forall (i,j) \in A_h, k \in K, c \in C \tag{23}$$

$$\widetilde{\delta}_{ijk}^c \in \{0,1\}, \quad \forall (i,j) \in A_r, k \in K, c \in C \tag{24}$$

$$Y_{sk}^c \in \{0,1\}, \quad \forall s \in S, k \in K, c \in C \tag{25}$$

The objective function (1) seeks to minimize the total system cost; specifically, the system cost includes the transportation cost on highway and railway links, the transfer cost at intermodal terminals, and the penalty cost for unsatisfied demands. Constraints (2) – (6) ensure flow conservation at highway nodes (*H*). Similarly, constraint (7) ensures flow conservation at railway nodes (*R*). Constraints (8) and (9) ensure flow conservation at intermodal terminals (*S*); constraint (8) maintains the conservation of flow if a terminal is selected whereas constraint (9) maintains the conservation of flow if the terminal is not selected. The decision variables $F_c^{sk}$ are calculated in constraint (10). Constraint (11) ensures that commodity shipments are delivered before the delivery deadline. Constraints (12) – (14) ensure that flows are less than or equal to the capacity of highway links, railway links, and intermodal terminals, respectively. Constraint (15) determines the unsatisfied demand. Lastly, constraints (16) – (18) are the relational constraints, constraints (19) – (21) are the definitional constraints, constraint (22) is the integrality constraint, and constraints (23) – (25) are the binary constraints. For constraints (16) – (18), the left-hand side term could be 0 instead of the product of $\varepsilon$. However, the formulation as presented provides a computational advantage. In addition, unsatisfied demands are assumed to be integer since the original demands are in number of intermodal containers.

As mentioned earlier, a transportation network may experience service disruptions. Hence, the above model (MIFR) with deterministic link capacities are not always valid. To account for uncertainty in the network, random capacity of highway link is denoted as $\hat{Q}_{ij}$, random capacity of railway link is denoted as $\hat{\widetilde{Q}}_{ij}$, and random capacity of intermodal terminal is denoted as $\hat{Q}_s$. Using these definitions, Eqs. (12), (13), and (14) have the following form.

$$\sum_{c \in C}\sum_{k \in K} d_k^c X_{ijk}^c \leq \hat{Q}_{ij}, \quad \sum_{c \in C}\sum_{k \in K} d_k^c \widetilde{X}_{ijk}^c \leq \hat{\widetilde{Q}}_{ij}, \quad \sum_{c \in C}\sum_{k \in K} d_k^c F_{sk}^c \leq \hat{Q}_s \tag{26}$$

To incorporate the modified constraints above into the optimization model, chance constraint programming is employed which guarantees that the solution satisfies the constraints over a subset of the sample space. Assume the following for a highway link capacity $\hat{Q}_{ij}$.

$$\hat{Q}_{ij} = Q_{ij}\left(1 + \lambda_{ij}\hat{\xi}_{ij}\right) \tag{27}$$

where $\lambda_{ij} \geq 0$ is a measure of uncertainty and $\hat{\xi}_{ij}$ represents a symmetric random variable on the interval [−1, 1]; meaning that $\hat{\xi}_{ij}$ and $-\hat{\xi}_{ij}$ have identical distributions. It should be noted that $\lambda_{ij}$ and $\hat{\xi}_{ij}$ is chosen in a way where $\hat{Q}_{ij} \geq 0$ always holds. Similarly, assume the following for the railway link and intermodal terminal capacities.

$$\hat{\widetilde{Q}}_{ij} = \widetilde{Q}_{ij}\left(1 + \widetilde{\lambda}_{ij}\hat{\widetilde{\xi}}_{ij}\right) \tag{28}$$

$$\hat{Q}_s = Q_s\left(1 + \lambda_s\hat{\xi}_s\right) \tag{29}$$

Let $E[Z]$ denotes the expected value of a random variable $Z$; then, $E[\hat{Q}_{ij}] = Q_{ij}$, $E\left[\hat{\widetilde{Q}}_{ij}\right] = \widetilde{Q}_{ij}$, and $E[\hat{Q}_s] = Q_s$. Hence, the model only requires the specification of mean values and the support of the random quantities instead of a specific probability distribution. Similar to chance constraint programming, this model has control over the likelihood that the constraints in Eq. (26) are violated. The following additional constraints are introduced in the model.

$$\sum_{c \in C}\sum_{k \in K} d_k^c X_{ijk}^c \leq Q_{ij} - \theta_{ij} \tag{30}$$

$$\sum_{c \in C}\sum_{k \in K} d_k^c \widetilde{X}_{ijk}^c \leq \widetilde{Q}_{ij} - \widetilde{\theta}_{ij} \tag{31}$$

$$\sum_{c \in C}\sum_{k \in K} d_k^c F_{sk}^c \leq Q_s - \theta_s \tag{32}$$

where $\theta_{ij} \geq 0$, $\tilde{\theta}_{ij} \geq 0$, and $\theta_s \geq 0$. The following probability expression

$$\Pr\left\{\sum_{c \in C}\sum_{k \in K} d_k^c X_{ijk}^c > \hat{Q}_{ij}\right\} \tag{33}$$

can be interpreted as the likelihood that the shipments based on the deterministic estimate of the highway link capacity exceed the realized capacity. To avoid this situation the probability in Eq. (33) needs to be acceptably small. Let us assume that $X_{ijk}^c$ is a feasible solution of the model defined by Eqs. (1)–(25), (30)–(32), then it follows

$$\Pr\left\{\sum_{c \in C}\sum_{k \in K} d_k^c X_{ijk}^c > \hat{Q}_{ij}\right\} = \Pr\left\{\sum_{c \in C}\sum_{k \in K} d_k^c X_{ijk}^c > Q_{ij}\left(1 + \lambda_{ij}\hat{\xi}_{ij}\right)\right\} < \Pr\left\{-Q_{ij}\lambda_{ij}\hat{\xi}_{ij} > \theta_{ij}\right\} \tag{34}$$

where the inequality follows from the following implication of events: since $X_{ijk}^c$ is feasible, then the event $\left\{\sum_{c \in C}\sum_{k \in K} d_k^c X_{ijk}^c > Q_{ij}\left(1 + \lambda_{ij}\hat{\xi}_{ij}\right)\right\}$ implies the event $\left\{-Q_{ij}\lambda_{ij}\hat{\xi}_{ij} > \theta_{ij}\right\}$. If the probability distribution of $\hat{\xi}_{ij}$ is assumed to be known, the right-most probability in Eq. (34) can easily be bounded. This is similar to a chance constraint programming approach.

Using a distribution-free approach (i.e., there is no explicit assumption about probability distributions), without loss of generality, let us assume that $\eta > 0$. Then it follows

$$\Pr\left\{-Q_{ij}\lambda_{ij}\hat{\xi}_{ij} > \theta_{ij}\right\} = \Pr\left\{\hat{\xi}_{ij} < \frac{\theta_{ij}}{-Q_{ij}\lambda_{ij}}\right\} = \Pr\left\{\hat{\xi}_{ij} > \frac{\theta_{ij}}{Q_{ij}\lambda_{ij}}\right\} = \Pr\left\{\exp(\eta\hat{\xi}_{ij}) > \exp\left(\frac{\eta\theta_{ij}}{Q_{ij}\lambda_{ij}}\right)\right\} \tag{35}$$

Markov's inequality gives the following equation from the last part of the above.

$$\Pr\left\{\exp(\eta\hat{\xi}_{ij}) > \exp\left(\frac{\eta\theta_{ij}}{Q_{ij}\lambda_{ij}}\right)\right\} \leq \exp\left(\frac{-\eta\theta_{ij}}{Q_{ij}\lambda_{ij}}\right) E\left[\exp(\eta\hat{\xi}_{ij})\right] \tag{36}$$

Since $\hat{\xi}_{ij}$ is a symmetric random variable, we can express $E\left[\exp(\hat{\eta}\xi_{ij})\right]$ as follows.

$$E\left[\exp(\hat{\eta}\xi_{ij})\right] = \int_{-1}^{0} \exp(\eta y)\, dF(y) + \int_{0}^{1} \exp(\eta y)\, dF(y) \tag{37}$$

$$= \int_0^1 \left[ \exp(\eta y) + \exp(-\eta y) \right] dF(y) \tag{38}$$

$$\leq \int_0^1 \max_{0 \leq y \leq 1} \left[ \exp(\eta y) + \exp(-\eta y) \right] dF(y) \tag{39}$$

$$\leq \left[ \exp(\eta) + \exp(-\eta) \right] \int_0^1 dF(y) \tag{40}$$

$$= \left[ \exp(\eta) + \exp(-\eta) \right] / 2 \tag{41}$$

Eq. (38) holds due to the symmetry of $\hat{\xi}_{ij}$ and Eq. (39) holds since the integrand is replaced by its maximum value. Inequality (40) is obtained by using the fact that integrand in Eq. (38) is maximized at $y = 1$. Again, using the symmetry Eq. (41) is obtained. Taylor series expansions of $\exp(\eta)$ and $\exp(-\eta)$ give us the following.

$$\left[ \exp(\eta) + \exp(-\eta) \right] / 2 \leq \exp(\eta^2 / 2) \tag{42}$$

Now, since $\eta > 0$ is arbitrary, the tightest possible bound can be obtained by minimizing over $\eta$. Therefore, using the above, Eq. (36) can be written as follows.

$$\Pr\left\{ \exp(\eta \hat{\xi}_{ij}) > \exp\left( \frac{\eta \theta_{ij}}{Q_{ij} \lambda_{ij}} \right) \right\} \leq \min_{\eta > 0} \exp\left( \frac{-\eta \theta_{ij}}{Q_{ij} \lambda_{ij}} \right) \exp(\eta^2 / 2) \tag{43}$$

To obtain Eq. (43) from the Eq. (36), it is assumed that random variations are symmetric, which is a common approach for solving robust optimization models (Bertsimas and Sim 2004; Ng and Waller 2012). The right-hand side of the Eq. (43) is strictly convex; hence, the unique optimal solution can be obtained by taking the derivative and setting it equal to zero. The optimal solution is

$$\eta^* = \frac{\theta_{ij}}{Q_{ij} \lambda_{ij}} \tag{44}$$

Substituting the above value into Eq. (43), the following can be obtained.

$$\Pr\left\{ \exp(\eta \hat{\xi}_{ij}) > \exp\left( \frac{\eta \theta_{ij}}{Q_{ij} \lambda_{ij}} \right) \right\} \leq \exp\left( \frac{-\theta_{ij}^2}{2(Q_{ij} \lambda_{ij})^2} \right) \tag{45}$$

The above discussion is summarized in the following proposition.

**Proposition 1.1** If $\hat{\xi}_{ij}$ is a symmetric random variable with support [-1, 1] and $\theta_{ij} = \sqrt{-2\log(q_{ij})}Q_{ij}\lambda_{ij}$, where $0 < q_{ij} \leq 1$, then

$$\Pr\left\{\sum_{c \in C}\sum_{k \in K} d_k^c X_{ijk}^c > \hat{Q}_{ij}\right\} \leq q_{ij}$$

Likewise, by imposing the constraints (31) and (32), the following two propositions can be shown.

**Proposition 1.2** If $\tilde{\hat{\xi}}_{ij}$ is a symmetric random variable with support [-1, 1] and $\tilde{\theta}_{ij} = \sqrt{-2\log(\tilde{q}_{ij})}\tilde{Q}_{ij}\tilde{\lambda}_{ij}$, where $0 < \tilde{q}_{ij} \leq 1$, then

$$\Pr\left\{\sum_{c \in C}\sum_{k \in K} d_k^c \tilde{X}_{ijk}^c > \tilde{\hat{Q}}_{ij}\right\} \leq \tilde{q}_{ij}$$

**Proposition 1.3** If $\hat{\xi}_s$ is a symmetric random variable with support [-1, 1] and $\theta_s = \sqrt{-2\log(q_s)}Q_s\lambda_s$, where $0 < q_s \leq 1$, then

$$\Pr\left\{\sum_{c \in C}\sum_{k \in K} d_k^c F_{sk}^c > \hat{Q}_s\right\} \leq q_s$$

## 3 Numerical Experiments

To demonstrate the applicability of the proposed modeling framework, an actual road-rail freight transport network shown in Figure 2 was used. It covers all of the states in the Gulf Coast, Southeastern and mid-Atlantic regions of the U.S.: Texas, Oklahoma, Louisiana, Alabama, Mississippi, Arkansas, Georgia, Florida, South Carolina, North Carolina, Tennessee, Kentucky, Virginia, Maryland, West Virginia, and Delaware. The network has a total of 682 links (U.S. interstates and major highways and Class I railroads) and 187 nodes, including 44 intermodal terminals. The Freight Analysis Zone (FAZ) centroids from the Freight Analysis Framework version 3 (FAF3) database (Federal Highway Administration 2016) were treated as actual origins and destinations of commodity shipments. There is a total of 48 centroids in the study region. Origin-Destination (OD) pairs were constructed from these 48 FAZ centroids, and demands are obtained from the FAF3 database. The demand data were filtered to include only those commodities typically transported via intermodal (Cambridge Systematics 2007), and demands were converted into containers using an average load of 40,000 lbs per container. It was assumed that all commodities need to be delivered within 7 days. The transport cost on highways and railways were estimated to be $1.67 per mile per shipment (Torrey and Murray 2014) and $0.60 per mile per shipment (Cambridge Systematics 1995), respectively. The transfer cost at intermodal terminals was estimated to be $70 per shipment (Winebrake et al. 2008). Using free-flow speeds, the travel times on highway and railway links were calculated. The number of potential paths between an origin and a destination could be large. For that reason, only those paths that were less than or equal to five times the corresponding minimum path length were considered in the path set for a specific OD pair. This approach is valid since the discarded paths would not have satisfied the delivery deadline constraint.

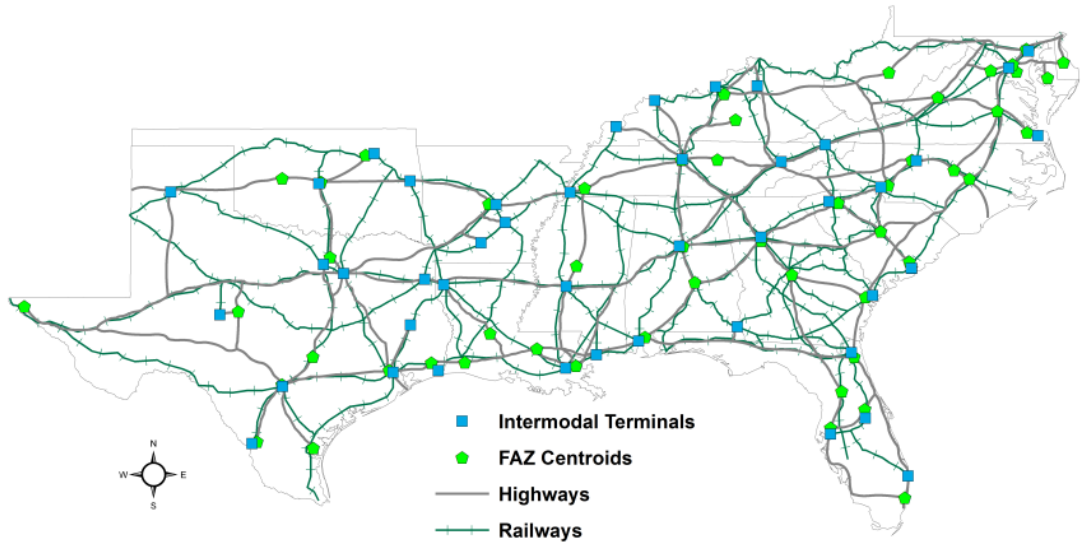

**Fig. 2** Large-scale U.S. road-rail intermodal network (Uddin and Huynh 2016)

To simulate network uncertainty, three types of disruptive events were considered in this study: (1) link disruption, (2) node disruption, and (3) intermodal terminal disruption. Note that affected links, nodes, or terminals are selected based on their vulnerability. A factorial experimental design (FED) was used to examine the effect of confidence level and capacity uncertainty parameters in the proposed model on total system cost (i.e., objective function value). In case of FED, "factors" and "levels" are utilized; "factors" are the variables that are chosen to be studied and "levels" are the predefined discrete values of the factors. The combination of all levels of factors are considered, and based on the resulting total system cost the effect of each combination of factors and levels is studied. Table 2 provides a summary of the FED. Three "factors" were considered: (1) number of disrupted elements, (2) confidence level $\left(q_{ij}, \tilde{q}_{ij}, q_s\right)$, and (3) capacity uncertainty $\left(\lambda_{ij}, \tilde{\lambda}_{ij}, \lambda_s\right)$. One can use the travel time reliability of highway links to infer its capacity uncertainty. For railway links and intermodal terminals, capacity uncertainty can be inferred from their disaster/disruption probabilities. Readers are referred to the works of Marufuzzaman et al. (2014) and Poudel et al. (2016) for more details regarding the methodology for estimating the disruption probability of intermodal network elements. For an experiment with a particular number of OD pairs and commodities, the combination of factors and levels result in a total of 112 instances for link disruptions, 112 instances for node disruptions, and 84 instances for intermodal terminal disruptions.

**Table 2** Summary of factorial experimental design

| Factors | Levels | | |
| --- | --- | --- | --- |
| | Link disruption | Node disruption | Terminal disruption |
| Number of disrupted elements | (1) 30, (2) 60, (3) 100, and (4) 200 | (1) 5, (2) 10, (3) 20, and (4) 40 | (1) 15, (2) 30, and (3) 44 |
| Confidence level $\left(q_{ij}, \tilde{q}_{ij}, q_s\right)$ | (1) 0.05, (2) 0.1, (3) 0.15, and (4) 0.2 | (1) 0.05, (2) 0.1, (3) 0.15, and (4) 0.2 | (1) 0.05, (2) 0.1, (3) 0.15, and (4) 0.2 |
| Capacity uncertainty $\left(\lambda_{ij}, \tilde{\lambda}_{ij}, \lambda_s\right)$ | (1) 0, (2) 0.05, (3) 0.1, (4) 0.15, (5) 0.2, (6) 0.25, and (7) 0.3 | (1) 0, (2) 0.05, (3) 0.1, (4) 0.15, (5) 0.2, (6) 0.25, and (7) 0.3 | (1) 0, (2) 0.05, (3) 0.1, (4) 0.15, (5) 0.2, (6) 0.25, and (7) 0.3 |

The proposed modeling framework was implemented in Python using Jetbrains PyCharm 5.0.3, and the IBM ILOG CPLEX 12.6 solver was used to solve the mixed-integer program. Experiments were run on a personal computer with Intel Core i7 3.20 GHz processor and 24.0 GB of RAM. For a given level of confidence and uncertainty level, using propositions 1.1–1.3, the amount of capacity reductions ($\theta$) can be obtained. Figures 3 to 6 present the experimental results for the real-world network for varying OD pairs and commodities.

Figure 3 depicts the resulting objective function values for 5 OD pairs (9 commodities) of shipments: (a) is for 30 disrupted links, (b) is for 60 disrupted links, (c) is for 100 disrupted links, (d) is for 200 disrupted links, (e) is for 5 disrupted nodes, (f) is for 10 disrupted nodes, (g) is for 20 disrupted nodes, (h) is for 40 disrupted nodes, (i) is for 15 disrupted intermodal terminals, (j) is for 30 terminals, and (k) is for 44 terminals. It can be seen that the objective function value increases with the level of uncertainty. Furthermore, increased confidence level leads to an increase in the objective function value. As expected, the objective function value increases as the number of affected links increases. A similar trend is observed for node and intermodal terminal disruptions. The objective function value was highest when all of the intermodal terminals were disrupted. This finding is logical because when all of the intermodal terminals are disrupted, commodities can only be shipped via road.

Figure 4 shows the variations of objective function values under different levels of capacity uncertainty and confidence levels for 10 OD pairs and 21 commodities. Figure 5 shows variations for 20 OD pairs and 43 commodities, and Figure 6 shows variations for 50 OD pairs and 87 commodities. The objective function values follow the same pattern observed in the 5 OD pairs scenario.

Collectively, the results indicate that under link and node disruption scenarios, most shipments are shipped via road-rail intermodal when a lower confidence level is considered. This can be attributed to the lower rail cost. When a higher confidence level is required under link and node disruptions, shipments are transported by road directly. This is can be attributed to the fact that a truck can always find an alternative route on the highway network when the intermodal network is disrupted. Freight shippers could use the above findings to make shipping decisions when the intermodal network is disrupted by some events. In summary, the managerial implications of the findings are that if freight shippers want a higher reliability for the delivery of its shipment under disruptions, they should ship via truck only. On the other hand, if reliability is not a concern, they should ship via road-rail intermodal due to lower cost.

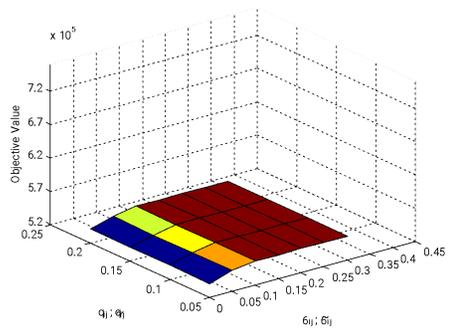 (a)
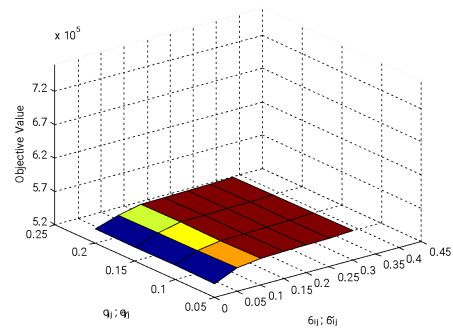 (b)
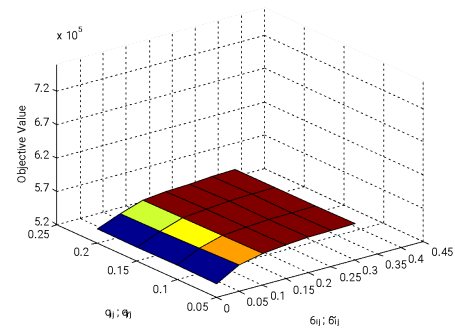 (c)
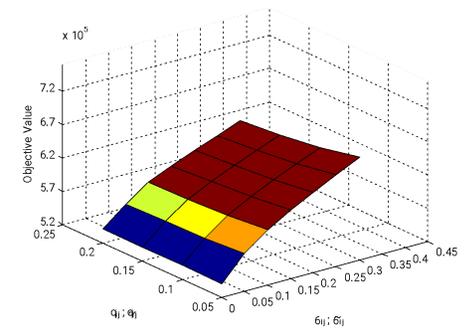 (d)
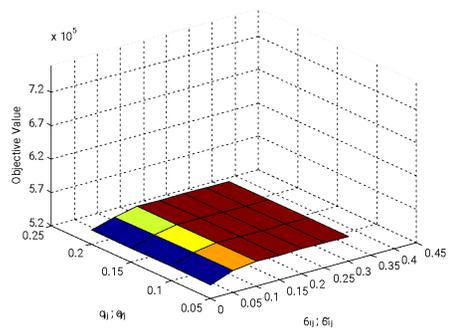 (e)
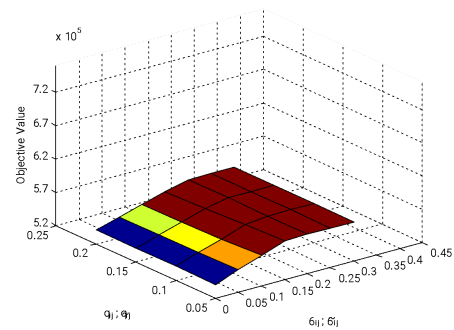 (f)
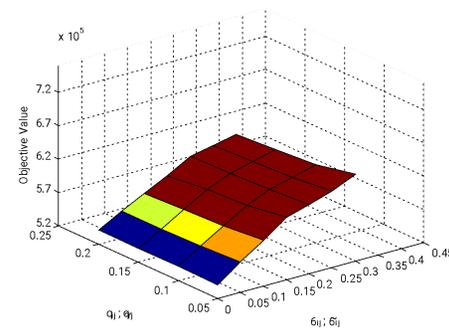 (g)
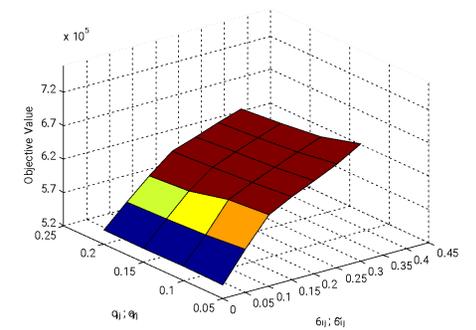 (h)
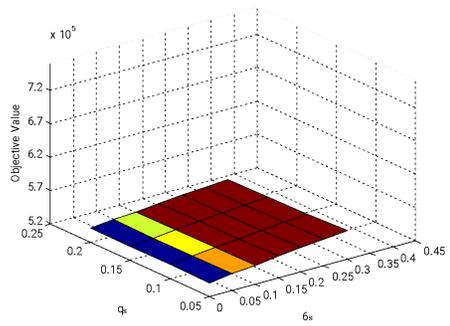 (i)
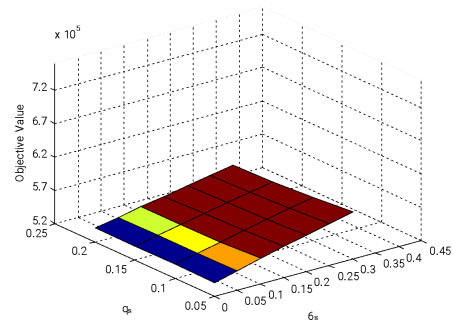 (j)
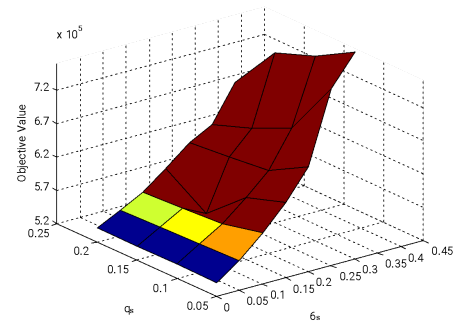 (k)

**Fig. 3** Objective function values under different levels of capacity uncertainty and confidence levels for 5 OD pairs (9 commodities): (a) 30 links, (b) 60 links, (c) 100 links, (d) 200 links, (e) 5 nodes, (f) 10 nodes, (g) 20 nodes, (h) 40 nodes, (i) 15 terminal, (j) 30 terminals, and (k) 44 terminals

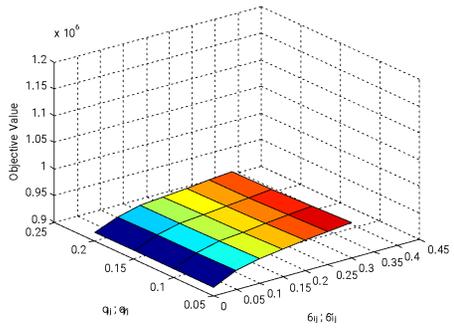
(a)

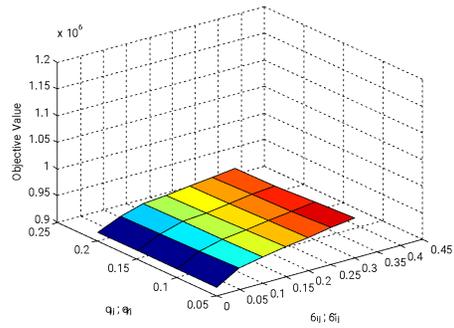
(b)

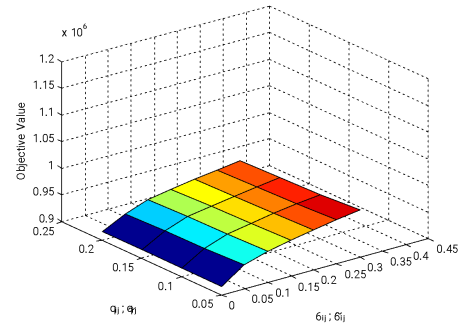
(c)

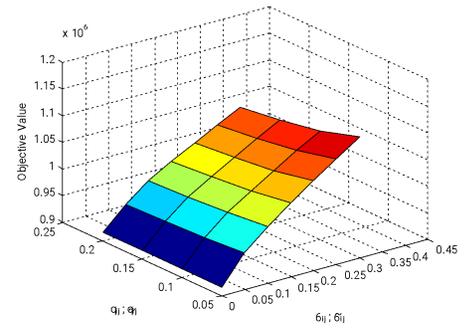
(d)

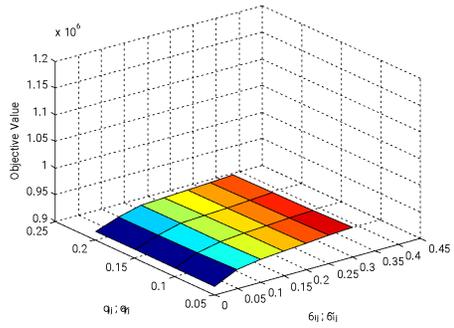
(e)

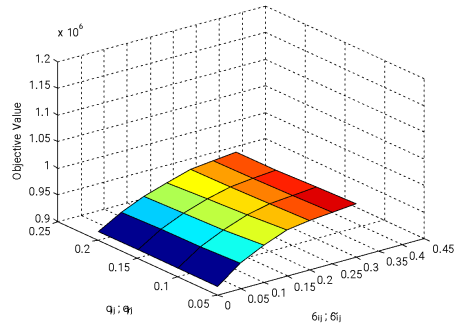
(f)

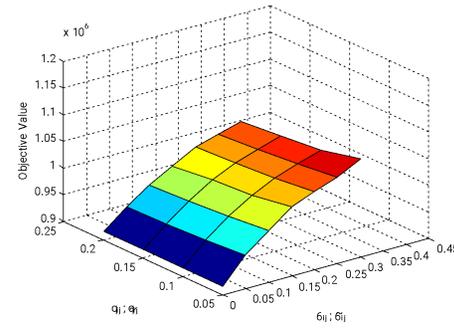
(g)

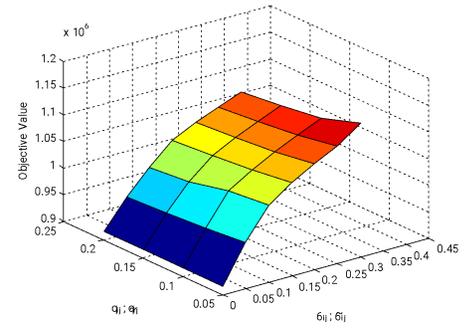
(h)

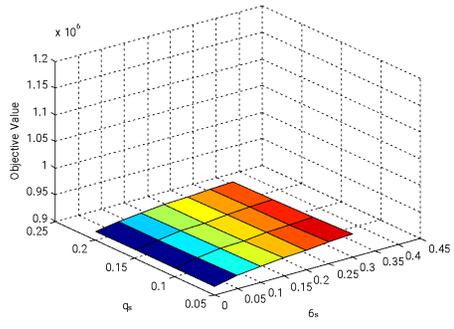
(i)

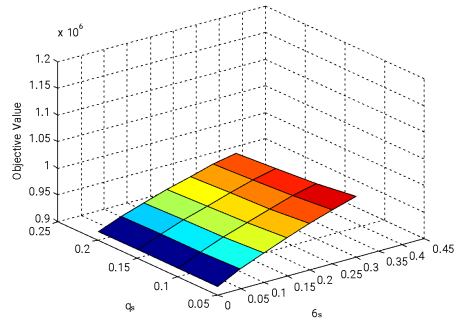
(j)

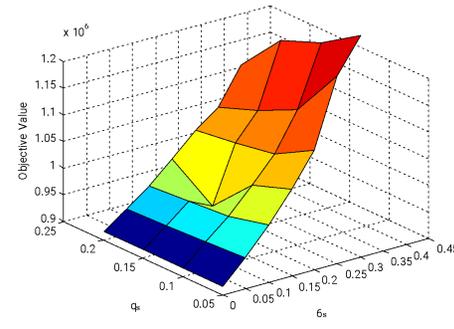
(k)

(i) (j) (k)

**Fig. 4** Objective function values under different levels of capacity uncertainty and confidence levels for 10 OD pairs (21 commodities): (a) 30 links, (b) 60 links, (c) 100 links, (d) 200 links, (e) 5 nodes, (f) 10 nodes, (g) 20 nodes, (h) 40 nodes, (i) 15 terminal, (j) 30 terminals, and (k) 44 terminals

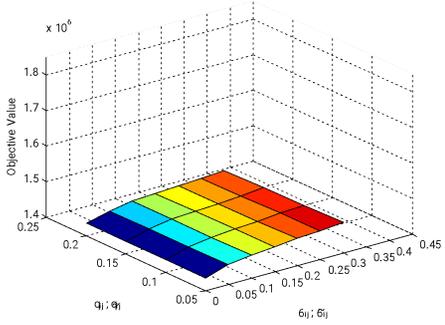
(a)

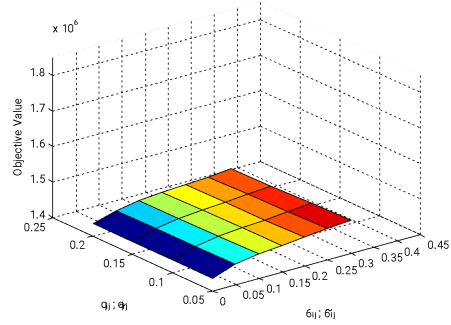
(b)

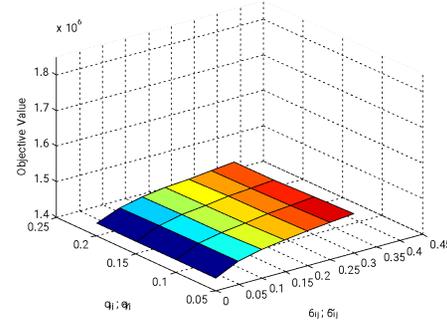
(c)

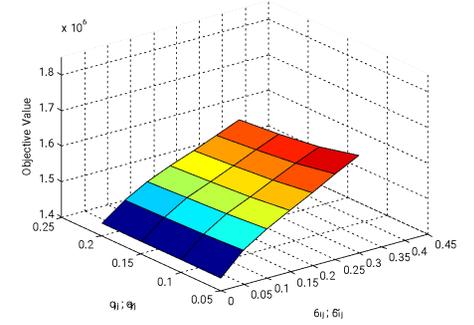
(d)

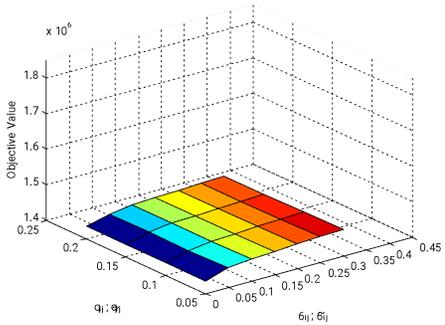
(e)

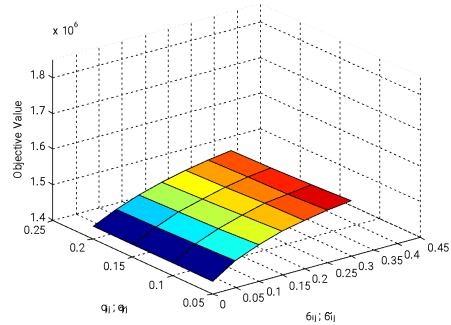
(f)

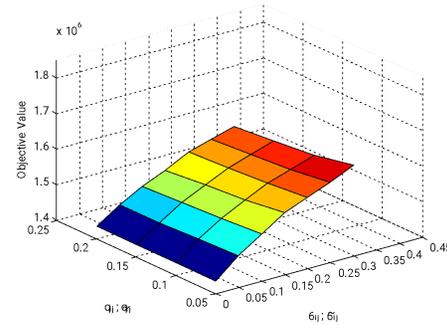
(g)

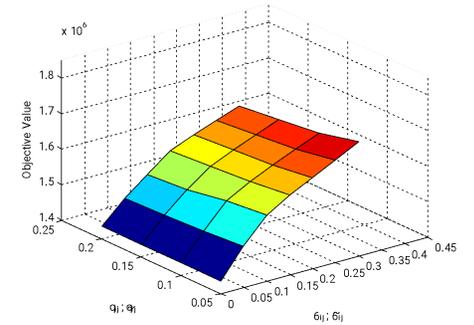
(h)

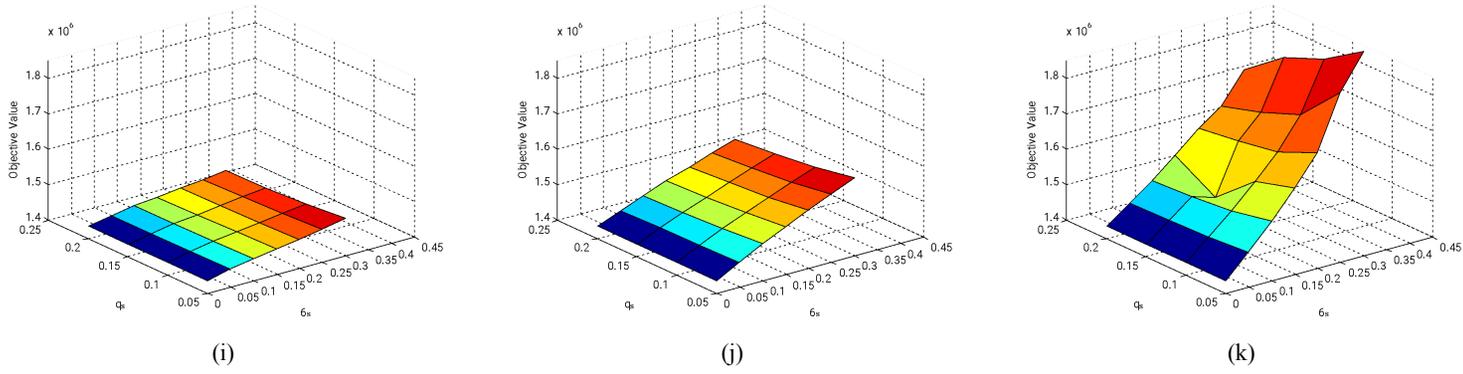

**Fig. 5** Objective function values under different levels of capacity uncertainty and confidence levels for 20 OD pairs (43 commodities): (a) 30 links, (b) 60 links, (c) 100 links, (d) 200 links, (e) 5 nodes, (f) 10 nodes, (g) 20 nodes, (h) 40 nodes, (i) 15 terminal, (j) 30 terminals, and (k) 44 terminals

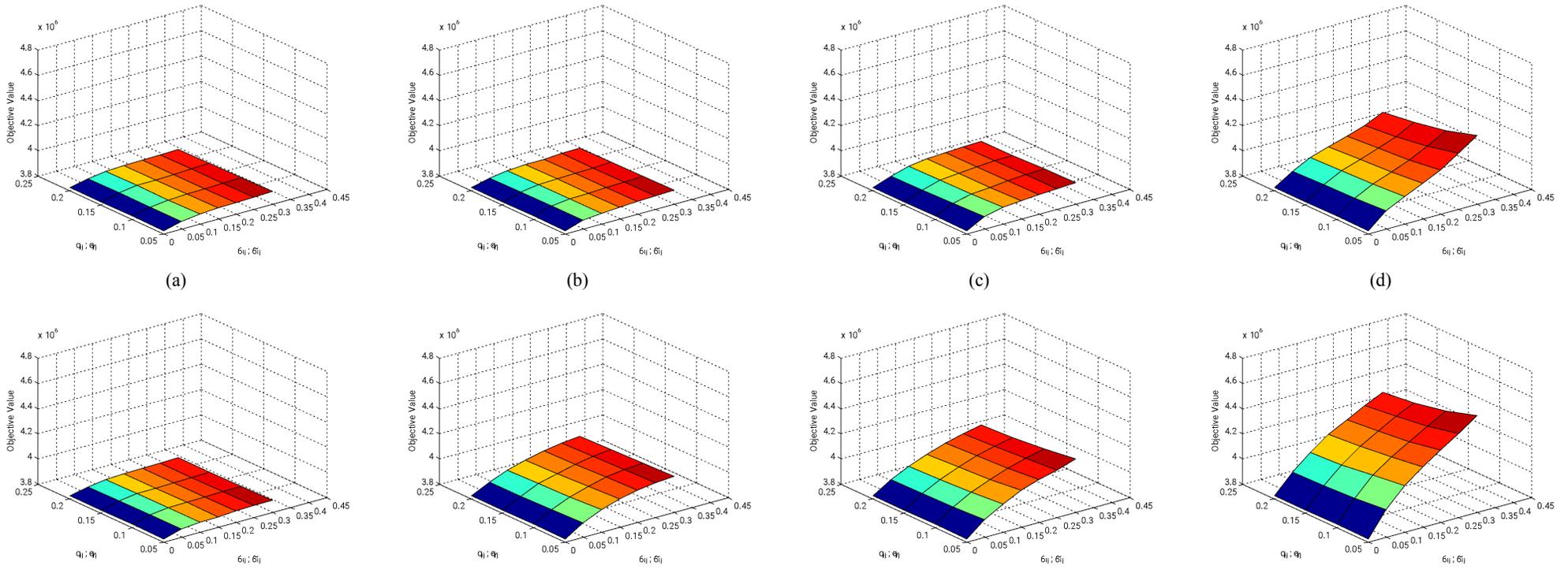

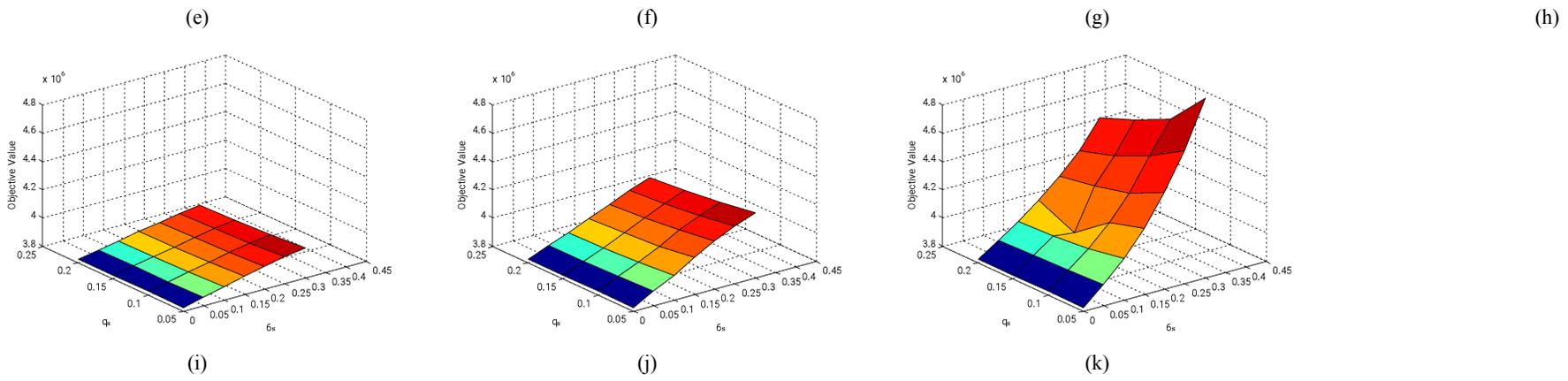

**Fig. 6** Objective function values under different levels of capacity uncertainty and confidence levels for 50 OD pairs (87 commodities): (a) 30 links, (b) 60 links, (c) 100 links, (d) 200 links, (e) 5 nodes, (f) 10 nodes, (g) 20 nodes, (h) 40 nodes, (i) 15 terminal, (j) 30 terminals, and (k) 44 terminals

After analyzing the experimental results, it is possible to quantify and compare vulnerability of different elements in road-rail intermodal freight transport networks. The observations are summarized in two propositions. Before presenting the propositions, an index called *importance* is introduced. It is assumed that a link is disrupted when its travel time ($\hat{t}_{ij}$) is greater than the typical travel time ($t_{ij}$). In particular, the term $\hat{t}_{ij}$ denotes the link travel time (with uncertainty) under the disruption scenario, and the term $t_{ij}$ denotes the link travel time (without uncertainty) under the normal scenario. The importance of a link $(i,j) \in A$ with respect to the entire network is defined as follows. The term $d_k^c X_{ijk}^c$ is a weight for each commodity $k \in K$ and OD pair $c \in C$ combination.

$$L_{ij}^{net} = \frac{\sum_{c \in C} \sum_{k \in K} d_k^c X_{ijk}^c (\hat{t}_{ij} - t_{ij})}{\sum_{c \in C} \sum_{k \in K} d_k^c X_{ijk}^c t_{ij}}, \qquad (i,j) \in A \tag{46}$$

In this study, vulnerability is defined in terms of reduced serviceability. It is possible to measure the reduced serviceability (i.e., vulnerability) by computing the increase in generalized cost of travel (i.e., travel time) for commodity shipments (Jenelius et al. 2006). To measure and compare the vulnerability of transportation network elements, a number of measures have been developed. These include *criticality* (Jenelius et al. 2006), *importance* (Jenelius et al. 2006; Rupi et al. 2015; Darayi et al. 2017), and *exposure* (Jenelius et al. 2006). The importance is defined above as the consequences of a network element in the road-rail intermodal network being disrupted. It is computed by accounting for the increase in travel time for each network element which in turn affects the performance of the network. For this reason, importance can be used to measure and compare the vulnerability of a network element. Furthermore, this index can be used to compare vulnerability across different types of network elements.

For a node disruption, all links connected to that node are affected. If set $J$ includes all the nodes connected to a specific node $i'$, i.e., $J = \{j \mid (i',j) \in A, i' \neq j\}$, then the importance of node $i'$ with respect to the entire network is defined as follows.

$$N_{i'}^{net} = \sum_{j \in J} L_{i'j}^{net} \tag{47}$$

The following two propositions apply to intermodal freight transport networks.

**Proposition 2.** *Impact of a node disruption is always greater than that of a link disruption if and only if both elements are affected by the same disruptive event.*

*Proof.* Suppose that impact of a network element disruption can be quantified by the importance measures defined above. Thus, the network element that has a higher importance value will have a greater impact on an intermodal freight transport network during a disruptive event. During a disruption, if the relative increase in link travel time ($\hat{t}_{ij}/t_{ij}$) is the same for all links inside the affected region or area, then by definition, $|J| \geq 2$, and hence, the following always holds for any node $i$ ($i = i'$ or $i \neq i'$).

$$N_{i'}^{net} > L_{ij}^{net}, \quad \forall (i,j) \in A \tag{48}$$

**Proposition 3.** *Impact of an intermodal terminal disruption is always greater than that of a node disruption if and only if both elements are affected by the same disruptive event.*

*Proof.* The importance of a terminal with respect to the entire network is defined as follows.

$$N_s^{net} = \sum_{j' \in J'} L_{sj'}^{net} + L_{sd}^{net} \tag{49}$$

The set $J'$ includes all the nodes connected to the terminal $s$ except for the dummy node $d$, i.e., $J' = J \setminus \{d\}$; to model a disruptive event, a dummy node and a dummy link is inserted between the terminal node and one of the network links connected to it. The travel time on the dummy link ($t_d$) is very large in the event of a disruption. Hence, the importance of any network node $i$ is always less than the importance of the terminal that it is connected to. Mathematically, this relationship can be expressed as follows.

$$N_s^{net} > N_i^{net}, \quad i \neq s \tag{50}$$

**Corollary 1.** *Impact of an intermodal terminal disruption is always greater than that of a link disruption if and only if both elements are affected by the same disruptive event.*

*Proof.* From Proposition 3, for any node $i'$ and intermodal terminal $s$, we have the following.

$$N_s^{net} > N_{i'}^{net} \tag{51}$$

From Proposition 2, for any node $i'$ and link $(i,j) \in A$, we have the following.

$$N_{i'}^{net} > L_{ij}^{net} \tag{52}$$

Hence, the following must always hold for any intermodal terminal $s$ and link $(i,j) \in A$.

$$N_s^{net} > L_{ij}^{net} \tag{53}$$

## 4 Summary and Conclusion

This paper proposed a new reliable modeling framework to determine the optimal routes for delivering multicommodity freight in an intermodal freight network that is subject to uncertainty. The finding from the proposed model is quite simple and intuitive: to ensure reliability, the model suggests that route planning be done by assuming the network elements have lower capacity than they actually have. To date, no formal framework has been developed to analytically determine the amount of capacity reduction needed to obtain a desired reliability level. This paper addressed this important gap by proposing a novel distribution-free approach. The framework is distribution-free in the sense that it only requires the specification of the mean values and the uncertainty intervals. The developed model was tested on an actual intermodal network in the Gulf Coast, Southeastern and mid-Atlantic regions of the U.S. It is found that the total system cost increases with the level of capacity uncertainty and with increased confidence levels for disruptions at links, nodes, and intermodal terminals.

This work can be improved in a couple of aspects. The link travel time was estimated based on its free-flow speed. A more accurate approach would be to use a link performance function; it should be noted that such an approach will lead to a non-linear model. The developed model considers only two modes: road and rail. A natural extension of this work would be to include other modes such as air and water. The impact of a disruption is limited to the disrupted source and its connected network elements and that impact of the disruption is assumed to be uniform. A more accurate approach would be to consider how a disruption propagates throughout the affected region and the different levels of impact on the affected network elements. Lastly, this work assumes that the random variable $\hat{\xi}_{ij}$ is symmetric. Future work could develop a model that does not depend on the symmetry of $\hat{\xi}_{ij}$.